\documentclass[showpacs,showkeys]{revtex4}
\usepackage[dvips]{graphicx}
\usepackage{epsfig}
\begin{document}

\title{Geodetic Coils on Deformed Sphere}

\author{V. L. Golo$^1$}
\email{golo@mech.math.msu.su }

\author{D. O. Sinitsyn$^1$}

\affiliation{$^1$ Department of Mechanics and Mathematics \\
     Moscow University \\
     Moscow 119 899   GSP-2, Russia \\  }

\date{July, 8, 2005}

\begin{abstract}
   We study geodetic lines on a
   surface generated by a small deformation of the standard
   2D-sphere.
   We construct an auxiliary hamiltonian system with the view of
   describing geodetic coils and almost closed geodesics, by using
   the fact that loops of the coil can be well approximated by
   great circles of the sphere.
   The phase space of  the auxiliary system is determined by
   the graph generated by separatrixes of its solutions,
   the vertices of the graph corresponding to
   almost closed geodesics and the edges to the geodetic coils
   joining them.  Topological types of the graph depend on the parameters
   determining the deformation. Using the method of averaging
   in conjunction with the computer  modelling of the auxiliary system,
   we obtain a fairly detailed visualization of geodesics on the deformed sphere.
\end{abstract}

\pacs{1111} \keywords{geodetic lines, averaging method,
separatrixe}

\maketitle

    \section{Geodetic coils}
    \label{introduction}

Geodesic lines on a surface can be considered either as straight
lines with respect to a Riemann metric, or trajectories of a
particle of mass $m$ moving freely on the surface. The second
approach allows for using the methods of analytical mechanics, and
has drawn considerable attention, \cite{arnold}. It should be
noted that the general solution of the problem of geodetic lines
is known only for certain special cases, e.g. ellipsoid,
\cite{arnold}, \cite{ndf}, and  the topology, or analysis situs,
of geodesics on a surface needs specific studying. In this paper
we use  asymptotic methods  and in particular focus our attention
on geodetic lines that are closed to within the first order of
perturbation theory.

\begin{figure}
  \begin{center}
    \includegraphics[width = 200bp]{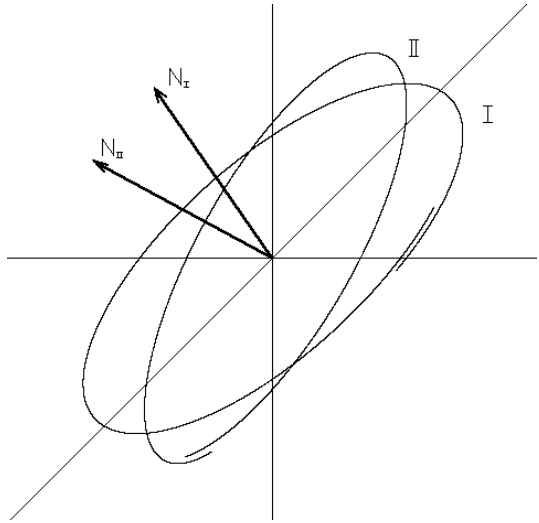}
    \caption{Two coils I, II of a geodesic on ellipsoid;
         vectors $ N_I $ and $ N_{II} $ are the normals to the planes of
         great circles approximating the coils}   \label{fig1}
  \end{center}
\end{figure}

\begin{figure}
  \begin{center}
    \includegraphics[width = 200bp]{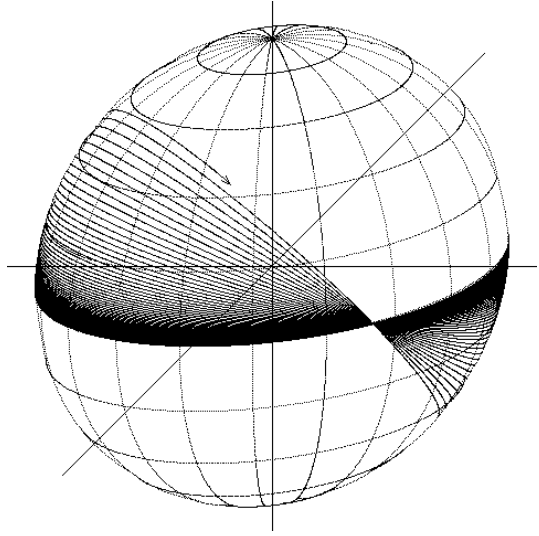}
    \caption{Ellipsoid with main axes 1.01, 1.02, 1.03;
         the geodesic, corresponding to a separatrixe of
         the auxiliary system, leaves a saddle point;
         the initial point (-0.0117, 0.0001, -1.0299);
         the initial velocity  (-1.9416, 0.0194, 0.0228.)
        }
    \label{fig2}
  \end{center}
\end{figure}

\begin{figure}
  \begin{center}
    \includegraphics[width = 200bp]{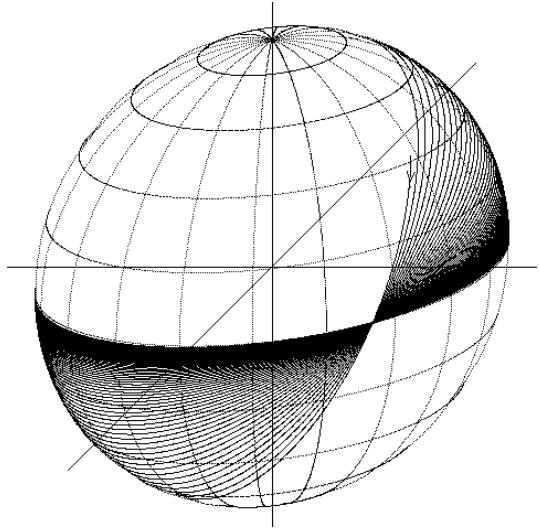}
    \caption{Ellipsoid with main axes 1.01, 1.02, 1.03;
             the geodesic, corresponding to a separatrixe of
             the auxiliary system, arrives at a saddle point;
             the initial point (0.7416, 0.6169, 0.3176);
             the initial velocity  (0.4728, -1.2782, 1.3832)
            }

    \label{fig3}
  \end{center}
\end{figure}

The central idea relies on the circumstance that in case of an
ellipsoid not differing substantially from a sphere, the great
circles of the latter may serve a good approximation to the
ellipsoid's geodesics, if they are short enough, and one can
visualize the geodesics as winding up in coils; loops or rings of
the coil corresponding to great circles of the sphere, see
FIG.\ref{fig1}. Hence approximating the successive rings by great
circles, we may describe the change in the position of the rings
by the motion of a great circle, see FIG.\ref{fig1}, which in its
turn is determined by the normal vector $\vec L$ of the plane
cutting the sphere along the great circle. To cast this picture in
a more quantitative form we may use the fact that the normal
vector $\vec L$ is the angular momentum of the particle moving
along the geodetic line. This, we obtain the advantage of allowing
for the use of asymptotic methods, in particular the averaging.
Our main instrument is an auxiliary hamiltonian system describing
the asymptotic motion of $\vec L$. It can be easily solved, so
that its stationary solutions correspond to closed planar
geodesics within the limit of precision provided by the averaging
method , see FIG.\ref{fig2}. We can get some insight into the
behaviour of geodesics by considering solutions joining the
stationary solutions, that is separaterixes, which, from a purely
geometrical point of view, are a kind of bridges between closed
geodesics. In the case of two dimensional torus in three
dimensional spcae this phenomenon had been  observed by
Yu.S.Volkov, \cite{volkov}. Generally, a  separaterixe geodetic
coil does not coincide anywhere with a closed geodesic but comes
infinitely close to it, and thus resembles a limit cycle, see
FIG.\ref{fig2}. The set of stationary points and separaterixes of
the auxiliary system forms a net, or graph, which, as is shown in
Section II can be realized on projective plane. Thus, we obtain a
topological classification of geodetic coils, and find that the
number of their topological types is finite.

These arguments  are valid even for surfaces having a more general
form than that of the ellipsoid. Studying the specific case of the
deformation determined by fourth order terms involves the surface
being substantially different from the ellipsoid as regards its
differential geometry, as well as the algebraic structure of its
equations, and hence we come across considerable difficulties by
attempting a direct investigation of the geodesics. The asymptotic
approach outlined above has the advantage of getting round this
problem.

    \section{Averaged equations of geodesics}
    \label{main_eq}

The equations determining geodesics on a surface can be cast in
the form of the equation, \cite{jacobi},

  \begin{equation}
    \ddot{\vec x}
    = \lambda \,
    \frac{\partial \varphi}
     {\partial \vec x}   \mbox{,}
    \label{2Newton}
  \end{equation}

\noindent where $ \varphi(\vec x) $ is the right-hand side of the
constraint $ \varphi(\vec x) = 0 $ determining the surface. The
Lagrangian multiplier can be found explicitly , so that the
equation of motion, in the form that does not involve $\lambda$,
reads

  \begin{equation}
    \ddot{\vec x} = - \frac{\dot{\vec x} \cdot
                  \displaystyle
                  \frac{\partial^2 \varphi}{\partial
                            \vec x^2} \cdot
                   \, \dot{\vec x}
                  }
             {\left( \displaystyle
                 \frac{\partial \varphi}{\partial \vec x}
              \right)^2} \,
	     \frac{\partial \varphi}{\partial \vec x}    \mbox{.}
    \label{2Newt_last}
  \end{equation}

\noindent In this paper we consider surfaces that do not differ
substantially from sphere, and in assume that their equations be
of the form

$$
  \varphi(\vec x) = \sum_i ( x_i^2  +  \varepsilon_i x_i^4 ) - 1
  \mbox{.}
$$

\noindent where $\varepsilon_i$ are small. Then
equations (\ref{2Newt_last}) read

  \begin{equation}
    \ddot{x_i} = - \frac{
     \sum_j (2 + 12 \varepsilon_j x_j^2) \dot{x_j}^2}
	{
     \sum_j (2 x_j + 4 \varepsilon_j x_j^3) ^2}
       (2 x_i + 4 \varepsilon_i x_i^3)   \mbox{.}
    \label{ConcrNewt}
  \end{equation}

The equations given above  are more tractable than the usual ones
employing the Christoffel symbols and an explicit parametrization
of sphere, so that one may prefer them for the needs of numerical
simulation, as is done in this paper.
With the view of
obtaining a qualitative description of geodesics, we shall
consider  the angular momentum
$$
   \vec L = \vec r \times \vec p,
$$

Its components satisfy the equations, which follow from
(\ref{ConcrNewt})

  \begin{equation}
    \begin{array}{lcl}
    \dot{L_1} = - 4 \, \displaystyle \frac{
     \sum_j (2 + 12 \varepsilon_j x_j^2) \dot{x_j}^2}
       {
     \sum_j (2 x_j + 4 \varepsilon_j x_j^3) ^2} \,
       x_2 x_3 (\varepsilon_3 x_3^2 - \varepsilon_2 x_2^2)
    \vspace{2mm} \\
    \dot{L_2} = - 4 \, \displaystyle \frac{
     \sum_j (2 + 12 \varepsilon_j x_j^2) \dot{x_j}^2}
       {
     \sum_j (2 x_j + 4 \varepsilon_j x_j^3) ^2} \,
       x_3 x_1 (\varepsilon_1 x_1^2 - \varepsilon_3 x_3^2)
    \vspace{2mm} \\
    \dot{L_3} = - 4 \, \displaystyle \frac{
     \sum_j (2 + 12 \varepsilon_j x_j^2) \dot{x_j}^2}
       {
     \sum_j (2 x_j + 4 \varepsilon_j x_j^3) ^2} \,
       x_1 x_2 (\varepsilon_2 x_2^2 - \varepsilon_1 x_1^2)
    \vspace{2mm}
    \end{array}
    \label{ConcrMom}
  \end{equation}

\noindent  Even though the equations given above are exact, in the
sense that they do not involve any approximation and do not use
the $\varepsilon_i$ being small, their treatment still need
further refining, and this will be done with the help of the
method of averaging. Generally, the  approach relies on studying
the evolution equations for integrals of motion of the unperturbed
system, i.e. in our case the normals to the planes of the large
circles, with respect to the basic periodic solution of the
latter. The averaging serves as a filter separating the main
regular part of the solution from the oscillating one caused by
small terms considered as perturbation, see \cite{hamming}.

We shall write the basic equation for the particle's  motion on
the sphere of unit radius in the form

$$
    \vec x = cos(\omega t + \theta) \vec e_1
       + cos(\omega t + \theta) \vec e_2
$$

\noindent vectors $\vec e_1, \vec e_2, \vec e_3$ determined by the
equations:

$$
\begin{array}{lcl}
  \vec e_1 = \displaystyle \frac{1}{\sqrt{L_2^2 + L_3^2}}(0, L_3, -L_2)
    \vspace{2mm} \\
  \vec e_2 = \displaystyle \frac{1}{L \sqrt{L_2^2 + L_3^2}}
    (-L_2^2 - L_3^2, L_1 L_2, L_1 L_3)
    \vspace{2mm} \\
  \vec e_3 = \displaystyle \frac{1}{L}(L_1, L_2, L_3) \mbox{.}
\end{array}
$$

The angular velocity $\omega$ is given by the
equation $\omega^2 = \dot{\vec x}^2 = L^2$, valid to within the
first order of perturbation.
Here $L_1, L_2, L_3$ are coordinates of the normal
to the plane of the great circle determining the solution, i.e.
the angular momentum.

Let us turn to the exact equations for the angular momentum
(\ref{ConcrMom}). With the help of the equations given above and
neglecting terms of the second, and higher, order in the $
\varepsilon_i $, we can transform equations (\ref{ConcrMom}) in
the form

  \begin{equation} 
    \begin{array}{rcl}
    \dot{L_1} = \displaystyle \frac{2 L^2 \varepsilon_2}{(L_2^2 + L_3^2)^2}
       \left[
	      \cos(\omega t + \theta) L_3
	      + \sin(\omega t + \theta) \displaystyle \frac{L_1 L_2}{L}
       \right]^3
       \left[
	       \cos(\omega t + \theta) (-L_2)
	       + \sin(\omega t + \theta) \displaystyle \frac{L_1 L_3}{L}
       \right] \vspace{2mm} \\
	       - \displaystyle \frac{2 L^2 \varepsilon_3}{(L_2^2 + L_3^2)^2}
       \left[
		\cos(\omega t + \theta) (-L_2)
		+ \sin(\omega t + \theta) \displaystyle \frac{L_1 L_3}{L}
       \right]^3
       \left[
		 \cos(\omega t + \theta) L_3
		  + \sin(\omega t + \theta) \displaystyle \frac{L_1 L_2}{L}
       \right] \\
    \end{array}
    \label{explicit}
  \end{equation}

 $$ 
    \begin{array}{rcl}
       \dot{L_2} = \displaystyle \frac{2 L^2 \varepsilon_3}{(L_2^2 + L_3^2)^2}
       \left[
	      \cos(\omega t + \theta) (-L_2)
	      + \sin(\omega t + \theta) \displaystyle \frac{L_1 L_3}{L}
       \right]^3
       \left[
	      \cos(\omega t + \theta) \cdot 0
	      + \sin(\omega t + \theta) \displaystyle \frac{- L_2^2 - L_3^2}{L}
       \right] \vspace{2mm} \\
	      - \displaystyle \frac{2 L^2 \varepsilon_1}{(L_2^2 + L_3^2)^2}
       \left[
	       \cos(\omega t + \theta) \cdot 0
	      + \sin(\omega t + \theta) \displaystyle \frac{- L_2^2 - L_3^2}{L}
       \right]^3
       \left[
		 \cos(\omega t + \theta) (-L_2)
		+ \sin(\omega t + \theta) \displaystyle \frac{L_1 L_3}{L}
       \right] \\
    \end{array}
 $$

 $$ 
    \begin{array}{rcl}
      \dot{L_3} = \displaystyle \frac{2 L^2 \varepsilon_1}{(L_2^2 + L_3^2)^2}
      \left[
	      \cos(\omega t + \theta) \cdot 0
	     + \sin(\omega t + \theta) \displaystyle \frac{- L_2^2 - L_3^2}{L}
      \right]^3
      \left[
		\cos(\omega t + \theta) L_3
	      + \sin(\omega t + \theta) \displaystyle \frac{L_1 L_2}{L}
      \right] \vspace{2mm} \\
	      - \displaystyle \frac{2 L^2 \varepsilon_2}{(L_2^2 + L_3^2)^2}
      \left[
		\cos(\omega t + \theta) L_3
	      + \sin(\omega t + \theta) \displaystyle \frac{L_1 L_2}{L}
     \right]^3
     \left[
		 \cos(\omega t + \theta) \cdot 0
	       + \sin(\omega t + \theta) \displaystyle \frac{- L_2^2 - L_3^2}{L}
      \right] \\
    \end{array}
  $$

It should be noted that the right-hand sides of Eqs.(\ref{explicit})
comprise terms oscillating in time and terms that vary slowly.
The situation can be treated within the framework of the averaging method,
\cite{hamming}, that is on neglecting the oscillatory terms we obtain the
averaged equations for the angular momentum

  \begin{equation}
    \begin{array}{rcl}
    \displaystyle
    \dot{L_1}
    &=& \displaystyle \frac34 \, \frac{L_2 L_3}{L^2}
    \left[ (\varepsilon_3 - \varepsilon_2) L_1^2 + \varepsilon_3 L_2^2
    - \varepsilon_2 L_3^2
    \right]    \mbox{,}
    \vspace{2mm}    \vspace{2mm} \\

    \displaystyle
    \dot{L_2}
    &=& \displaystyle \frac34 \, \frac{L_3 L_1}{L^2}
    \left[ - \varepsilon_3 L_1^2 + (\varepsilon_1 - \varepsilon_3) L_2^2
    + \varepsilon_1 L_3^2
    \right]    \mbox{,}
    \vspace{2mm}    \vspace{2mm} \\

    \displaystyle
    \dot{L_3}
    &=& \displaystyle \frac34 \, \frac{L_1 L_2}{L^2}
    \left[\varepsilon_2 L_1^2 - \varepsilon_1 L_2^2 +
    (\varepsilon_2 - \varepsilon_1) L_3^2
    \right]    \mbox{.}
    \vspace{2mm}            \\
    \end{array}
    \label{avmom}
  \end{equation}

It is worth noting that Eq.(\ref{avmom}) have
the Hamiltonian form determined by the usual Poisson brackets for
the angular momentum, \cite{routh}, \cite{arnold},
$$
    \{L_i, L_j\} = \sum_k \varepsilon_{ijk} L_k  \mbox{,}
$$
and the Hamiltonian

  \begin{equation}
    H = \frac{3}{16} L^2 \sum_i \varepsilon_i
             \left[ \left( \frac{L_i}{L} \right)^2 - 1
             \right]^2  \mbox{.}
    \label{hamiltonian}
  \end{equation}

\noindent This circumstance is particularly interesting because,
usually, the averaging procedure is not compatible with
Hamiltonian structure.

We may infer from  the two integrals of motion, $L^2$ and
$H$, that they admit of an explicit exact solution
that can be cast in the form of the equation

\begin{equation}
  t = \mp \frac43 \, (\varepsilon_2 + \varepsilon_3) \, L^2
  \int \frac{dL_1}
  {\sqrt{D (\varepsilon_2 L^2 - \varepsilon_3 L_1^2 \mp \sqrt{D})
       (\varepsilon_3 L^2 - \varepsilon_2 L_1^2 \pm \sqrt{D}) }}
  \mbox{.}
  \label{time}
\end{equation}
in which
$$
  D = - \varepsilon_1 (\varepsilon_2 + \varepsilon_3)(L^2 - L_1^2)^2
  - \varepsilon_2 \varepsilon_3 (L^2 + L_1^2)^2
  + \frac{16}{3} (\varepsilon_2 + \varepsilon_3) L^2 H
  \mbox{.}
$$

The exact solution is not particularly practical.
In this paper we shall not use it, but employ topological
considerations and numerical simulation so as to understand the
dynamics of $\vec L$ and, consequently, the behaviour of geodetic
coils.  The important point is considering the stationary
solutions to Eqs.(\ref{avmom}) for which the right-hand sides turn
out to be zero, and split into three parts {\sf S1, S2} and {\sf
S3}, determined by conditions on $\varepsilon_i$, as follows.

\begin{itemize}
  \item[\sf S1] \label{S1}
    No algebraic constraints imposed on $\varepsilon_i$:
    \begin{itemize}
        \item[\sf a.] \label{T1a}
        $ L_{10}  =  0, \quad L_{20}  =  0, \quad L_{30} \ne 0; $
        \item[\sf b.] \label{T1b}
        $ L_{10}  =  0, \quad L_{20} \ne 0, \quad L_{30}  =  0; $
        \item[\sf c.] \label{T1c}
         $ L_{10} \ne 0, \quad L_{20}  =  0, \quad L_{30}  =  0;$
    \end{itemize}
  \item[\sf S2] \label{S2}
    The constraints on $\vec L$ relaxed and linear constraints
    imposed on $\varepsilon_i$:
    \begin{itemize}
        \item[\sf a.] \label{T2a}
        $ L_{10} = 0, \quad L_{20} \ne 0, \quad L_{30} \ne 0, \quad
          \varepsilon_3 L_{20}^2 - \varepsilon_2 L_{30}^2 = 0  $
        \item[\sf b.] \label{T2b}
          $ L_{20} = 0, \quad L_{30} \ne 0, \quad L_{10} \ne 0, \quad
           \varepsilon_1 L_{30}^2 - \varepsilon_3 L_{10}^2 = 0  $
        \item[\sf c.] \label{T2c}
         $L_{30} = 0, \quad L_{10} \ne 0, \quad L_{20} \ne 0, \quad
         \varepsilon_2 L_{10}^2 - \varepsilon_1 L_{20}^2 = 0 $
    \end{itemize}
  \item[\sf S3] \label{S3}
    Vector $\vec L$ subject to  $L_{10} \ne 0, L_{20} \ne 0,  L_{30} \ne 0$
    and the quadratic constraints imposed on
    $\varepsilon_i$:
        $$
            \frac{L_{10}^2}
            {\varepsilon_1 \,\varepsilon_2 - \varepsilon_2 \,
            \varepsilon_3 + \varepsilon_3
            \,\varepsilon_1} =
            \frac{L_{10}^2}
            {\varepsilon_1 \,\varepsilon_2 + \varepsilon_2
            \,\varepsilon_3 - \varepsilon_3
            \,\varepsilon_1} =
            \frac{L_{10}^2}
            { - \varepsilon_1 \,\varepsilon_2 + \varepsilon_2
            \,\varepsilon_3 + \varepsilon_3
            \,\varepsilon_1}
        $$
 \end{itemize}

It is worth noting that equations  {\sf S2} involve the fulfilment
of the inequalities $\varepsilon_2 \varepsilon_3 > 0$,
$\varepsilon_3 \varepsilon_1 > 0$, and $\varepsilon_1
\varepsilon_2 > 0$ for cases {S2.a, S2.b, S2.c}, respectively,
whereas equations {\sf S3} involve
        \begin{eqnarray}
            \varepsilon_1 \,\varepsilon_2 - \varepsilon_2 \,\varepsilon_3 +
            \varepsilon_3
            \,\varepsilon_1  \, > 0  \nonumber  \\
            \varepsilon_1 \,\varepsilon_2 + \varepsilon_2 \,\varepsilon_3 -
            \varepsilon_3
            \,\varepsilon_1  \, > 0  \nonumber  \\
            \ - \varepsilon_1 \,\varepsilon_2 + \varepsilon_2 \,\varepsilon_3
            + \varepsilon_3
            \,\varepsilon_1  \, > 0  \nonumber
        \end{eqnarray}

Linearizing Eqs.(\ref{avmom}) at the stationary solutions and,
considering small fluctuations of $\vec L$ round them, we may
study their stability, which turns out to be determined by the
requirements
\begin{itemize}
    \item[\sf S1]
    \begin{itemize}
        \item[\sf a.] $\varepsilon_1 \varepsilon_2 > 0$;
        \item[\sf b.] $\varepsilon_2 \varepsilon_3 > 0$;
        \item[\sf c.] $\varepsilon_3 \varepsilon_1 > 0$.
    \end{itemize}
    \item[\sf S2]
    \begin{itemize}
        \item[\sf a.] $ \varepsilon_1 \varepsilon_2
            -\varepsilon_2 \varepsilon_3
            +\varepsilon_3 \varepsilon_1 < 0 $;
        \item[\sf b.] $ \varepsilon_1 \varepsilon_2
            +\varepsilon_2  \varepsilon_3
            -\varepsilon_3 \varepsilon_1 < 0 $;
        \item[\sf c.] $-\varepsilon_1 \varepsilon_2
            +\varepsilon_2 \varepsilon_3
            +\varepsilon_3 \varepsilon_1 < 0 $.
    \end{itemize}
    \item[\sf S3] \quad any $ \varepsilon_i $.
\end{itemize}

We may  put these equations in a more graphic form using the
integral $ L^2 = const $, and consider the motion of $\vec L$ on a
sphere of  fixed radius, the integral of energy $H$ taking
appropriate values. Then the stable solutions are fixed points
as regards Eqs.(\ref{avmom}),  the stable and the unstable points
are foci and saddle points, respectfully,
the separaterixes being lines joining the fixed points.
Together, they generate a graph on the sphere,
having the fixed points as vertices and
the separaterixes as edges.  It is important that the separaterixes,
i.e. the edges of the graph, are oriented according to the time,
$t$, so that the graph is the oriented one, and invariant with respect
to the symmetry $ \vec R \rightarrow - \vec R $, $ t \rightarrow - t $.

\begin{figure}
  \begin{center}
    \includegraphics[width = 200bp]{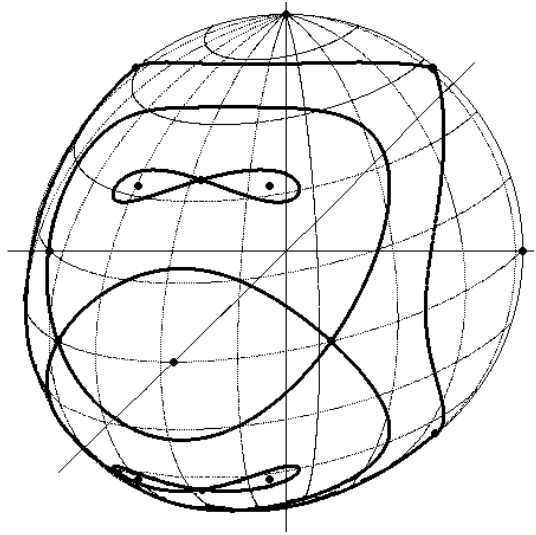}
    \caption{Separatrix net corresponding to the graph of
             {\sf Type I}  on the sphere; the twin points
             correspond to the symmetry given by Eq.(\ref{symmetry})
            }
            \label{fig4}
  \end{center}
\end{figure}

\begin{figure}
  \begin{center}
    \includegraphics[width = 200bp]{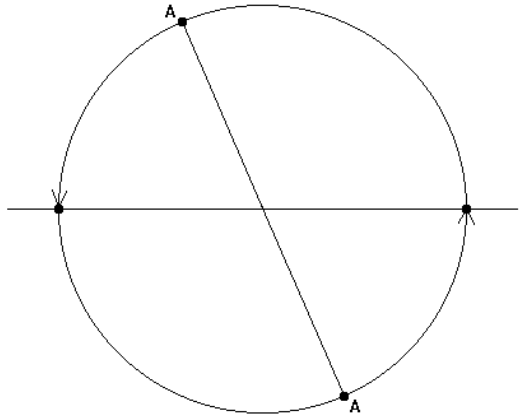}
    \caption{Diagramme of the projective plain by means of a disk
             with identified antipodal points at the boundary circle.
            }
            \label{fig5}
  \end{center}
\end{figure}

\begin{figure}
  \begin{center}
    \includegraphics[width = 200bp]{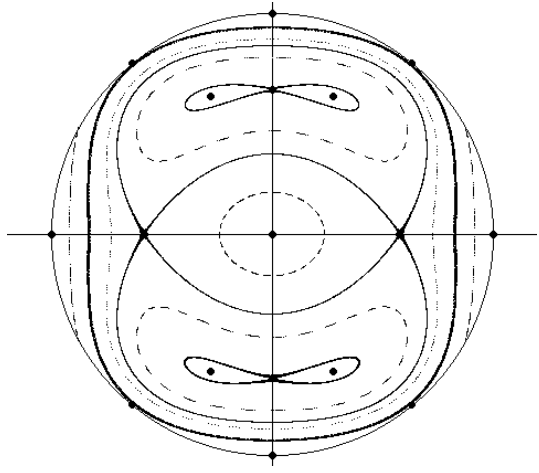}
    \caption{ {\sf Type I} trajectories of the auxiliary system on projective
          plain; dashed lines are typical trajectories,
          the solid ones separatrixes.
          Parameters of deformation:
          $ \varepsilon_1 = 0.02, \
        \varepsilon_2 = 0.03, \
        \varepsilon_3 = 0.04. $
        }
        \label{fig6}
  \end{center}
\end{figure}

\begin{figure}
  \begin{center}
    \includegraphics[width = 200bp]{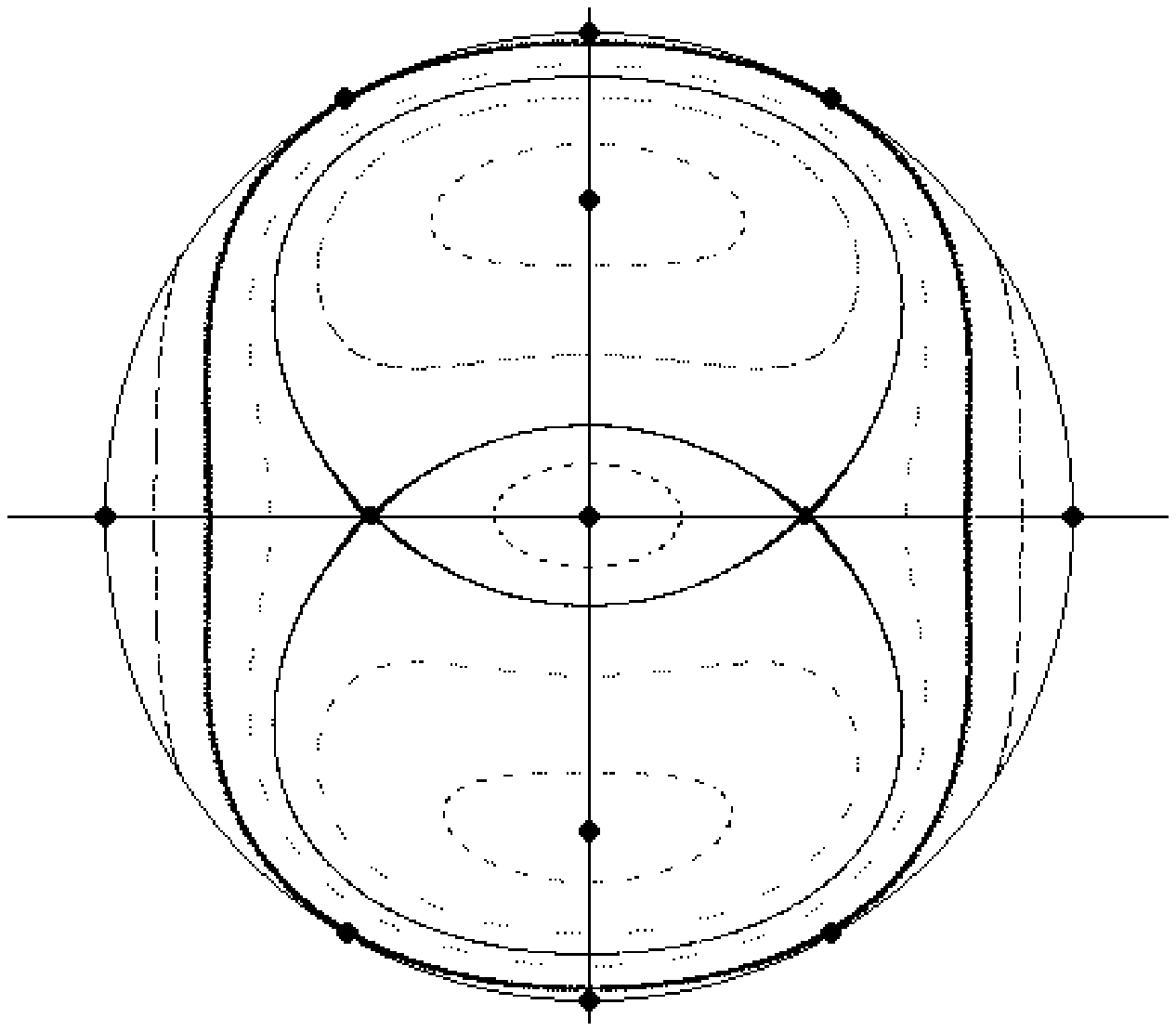}
    \caption{{\sf Type II} trajectories of the auxiliary system on the projective
          plain; dashed lines are typical trajectories,
          the solid ones separatrixes.
          Parameters of deformation:
          $ \varepsilon_1 = 0.01, \
        \varepsilon_2 = 0.03, \
        \varepsilon_3 = 0.04. $
        }
        \label{fig7}
  \end{center}
\end{figure}

The pictures of the separatrixe net on the sphere are rather
difficult for visualizing, and therefore there is a need for some
means which are easier to implement. We shall employ the familiar
construction of the projective plane, which runs as follows. Given
a pair of twin points $\vec R$ and $ - \vec R$ belonging to the
sphere we shall choose the representative of the pair determined
by the condition $R_3 \ge 0$; this enables us to consider only the
upper part of the sphere. Next take the projection of the upper
hemi-sphere on x-y plane along axes-z. Thus, a pair of antipodal
twins obtains a representative in the disk, i.e. a point inside
the disk or a pair of antipodal points at the boundary. The
separatrize nets constructed in this way on the projective plain,
are shown in FIGs(\ref{fig4} --- \ref{fig7}). It is important that
the equations of motion are invariant under the transformations
$P_i, \quad i=1,2,3$
\begin{equation}
   P_i : \quad t \rightarrow - t, \qquad  L_i \rightarrow -  L_i
   \label{symmetry}
\end{equation}
which generate automorphisms of the graph.

\noindent From the fact that the transformation
\begin{equation}
   P:    \quad t \rightarrow - t, \qquad  L_i \rightarrow -  L_i, \qquad
   i=1,2,3, \qquad P = P_1 P_2 P_3
   \label{sym2}
\end{equation}
leaves Eqs.(\ref{avmom}) invariant, we infer that a point
belonging to a solution of Eqs.(\ref{avmom}) has its counterpart,
or twin, at the antipodal point and therefore one may visualize
the dynamics of $\vec L$ on a sphere with identified antipodal
points, that is the projective plane.

Now we are in a position to determine the topological types of the
graphs by employing the computer simulation of the
Eqs.(\ref{avmom}) in conjunction with the knowledge of the
topological types of the fixed points obtained above.
It should be noted that we must check as to whether the
solutions provided by Eqs.(\ref{avmom}) agree with
those given by original Eqs.(\ref{2Newton}), see
FIG.(\ref{fig10}). The phase picture  can be obtained by
constructing a mesh generated by solutions to Eqs.(\ref{avmom}),
taking into account the types of fixed points. The results are
illustrated in FIGs(4 --- 7).

\begin{figure}
  \begin{center}
    \includegraphics[width = 200bp]{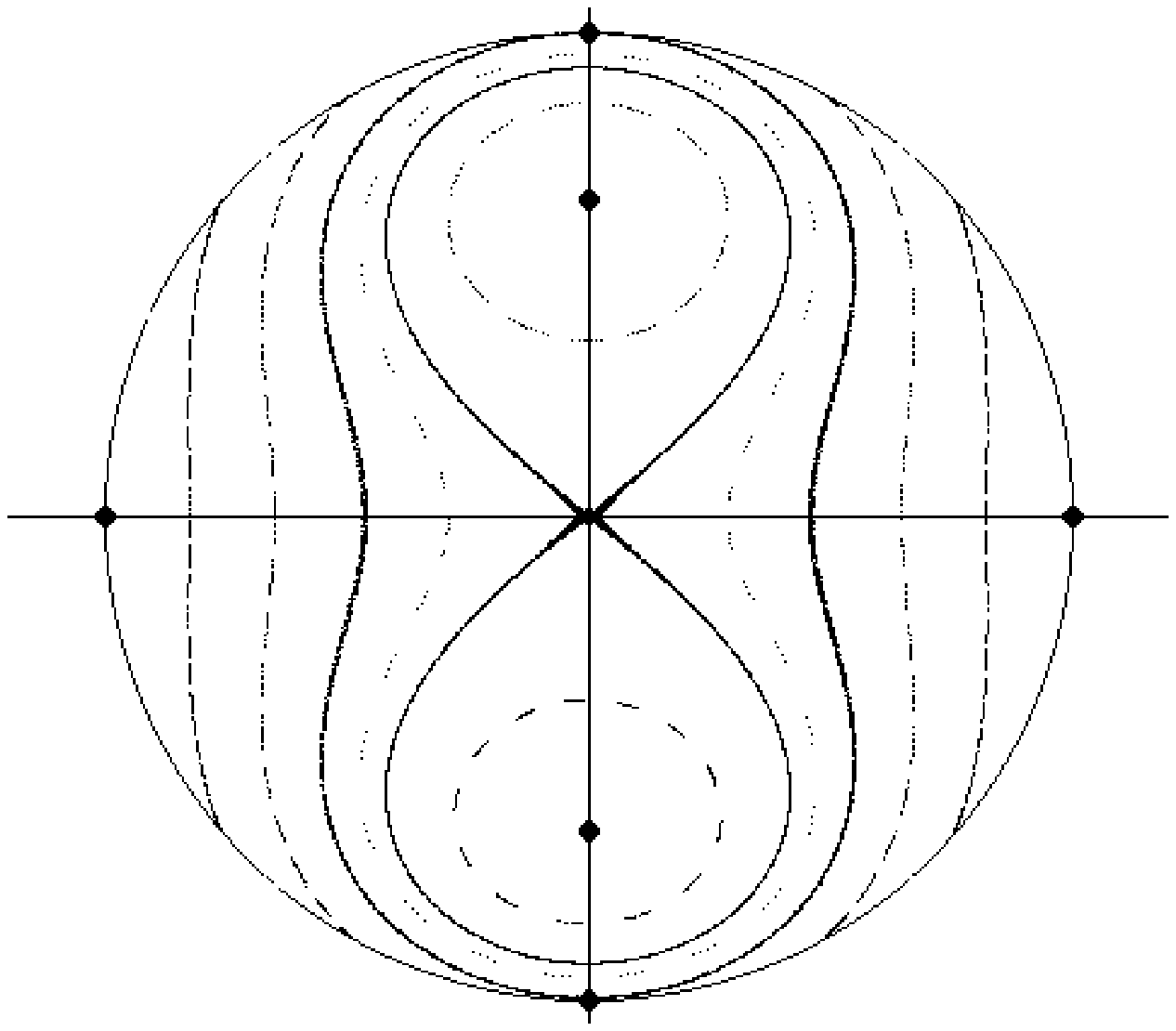}
    \caption{{\sf Type III} trajectories of the auxiliary system on the projective
          plain; dashed lines are typical trajectories,
          the solid ones separatrixes.
          Parameters of deformation:
          $ \varepsilon_1 = -0.02, \
        \varepsilon_2 = 0.03, \
        \varepsilon_3 = 0.04. $
        .}    \label{fig8}
  \end{center}
\end{figure}

\begin{figure}
  \begin{center}
    \includegraphics[width = 200bp]{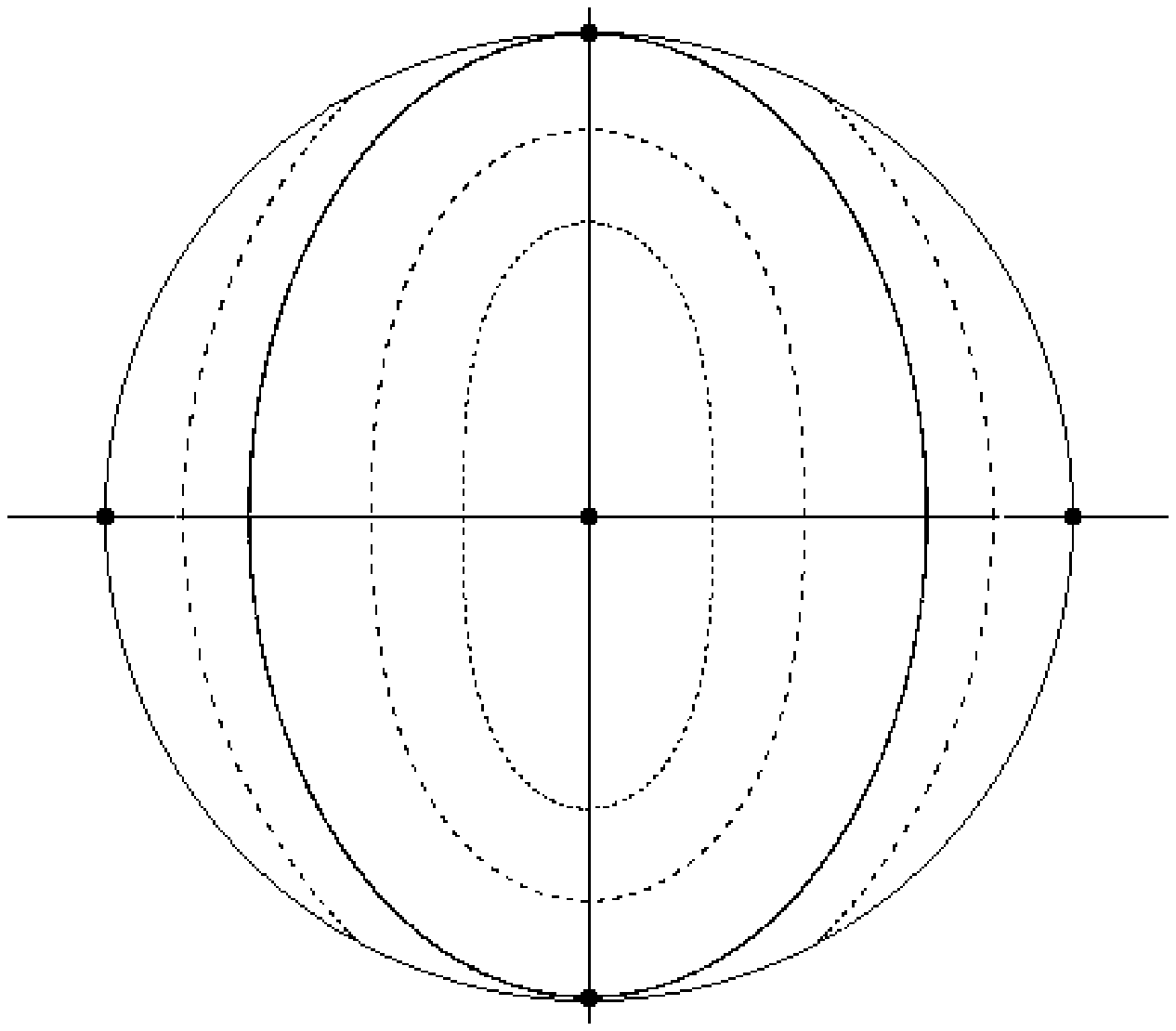}
    \caption{{\sf Type IV} trajectories of the auxiliary system on the projective
          plain; dashed lines are typical trajectories,
          the solid ones separatrixes.
          Parameters of deformation:
          $ \varepsilon_1 = -0.01, \
        \varepsilon_2 = 0.00, \
        \varepsilon_3 = 0.01. $
         }
         \label{fig9}
 \end{center}
\end{figure}

\begin{figure}
  \begin{center}
    \includegraphics[width = 200bp]{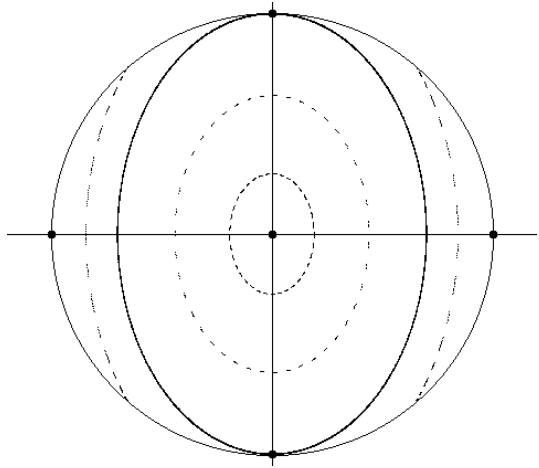}
    \caption{Ellipsoid with the main axes: 1.01, 1.02, 1.03;
         dashed lines are typical trajectories,
         the solid separatrixes.
         }
         \label{fig10}
  \end{center}
\end{figure}

We obtain the following topological types of the graphs:

\begin{itemize}
  \item[\sf Type I] \label{T1}
                    \quad FIG.\ref{fig4},
                    \quad  7 foci and 6 saddles;
                    $ \varepsilon_i $ being subject to the constraints:
      \begin{equation}
      \begin{array}{rcl}
         \varepsilon_1 \varepsilon_2
        -\varepsilon_2 \varepsilon_3
        +\varepsilon_3 \varepsilon_1 > 0 \\

     \varepsilon_1 \varepsilon_2
    +\varepsilon_2 \varepsilon_3
    -\varepsilon_3 \varepsilon_1 > 0 \\

    -\varepsilon_1 \varepsilon_2
    +\varepsilon_2 \varepsilon_3
    +\varepsilon_3 \varepsilon_1 > 0
      \end{array}
      \mbox{.}
     \label{t1constraints}
      \end{equation}

 \item[\sf Type II] \label{T2}
                   \quad FIG.\ref{fig5},
                   \quad 5 foci and 4 saddle points;
                   $ \varepsilon_i $ are not equal to zero,
                   have the same sign,
                   and at least one of
                   eqs.(\ref{t1constraints}) is not true.

 \item[\sf Type III] \label{T3}
                   \quad FIG.\ref{fig6},
                   \quad 3 foci and 2 saddle points;
                   $ \varepsilon_i $ being subject to
                   one of the following constraints:
           $ \varepsilon_2 \varepsilon_3 > 0 $ and
           $ \varepsilon_1 \varepsilon_2 \leq 0 $; \quad
           $ \varepsilon_3 \varepsilon_1 > 0 $ and
           $ \varepsilon_2 \varepsilon_3 \leq 0 $; \quad
           $ \varepsilon_1 \varepsilon_2 > 0  $ and
           $ \varepsilon_3 \varepsilon_1 \leq 0 $.

 \item[\sf Type IV] \label{T4}
                  \quad FIG.\ref{fig7},
                  \quad 2 foci and 1 saddle point;
                  $ \varepsilon_i $ being subject to
                  one of the following constraints:
           $ \varepsilon_1 = 0 $ and
           $ \varepsilon_2 \varepsilon_3 \leq 0 $; \quad
           $ \varepsilon_2 = 0 $ and
           $ \varepsilon_3 \varepsilon_1 \leq 0 $; \quad
           $ \varepsilon_3 = 0 $ and
           $ \varepsilon_1 \varepsilon_2 \leq 0 $.

  \end{itemize}

\begin{figure}
  \begin{center}
	\includegraphics[width = 400bp]{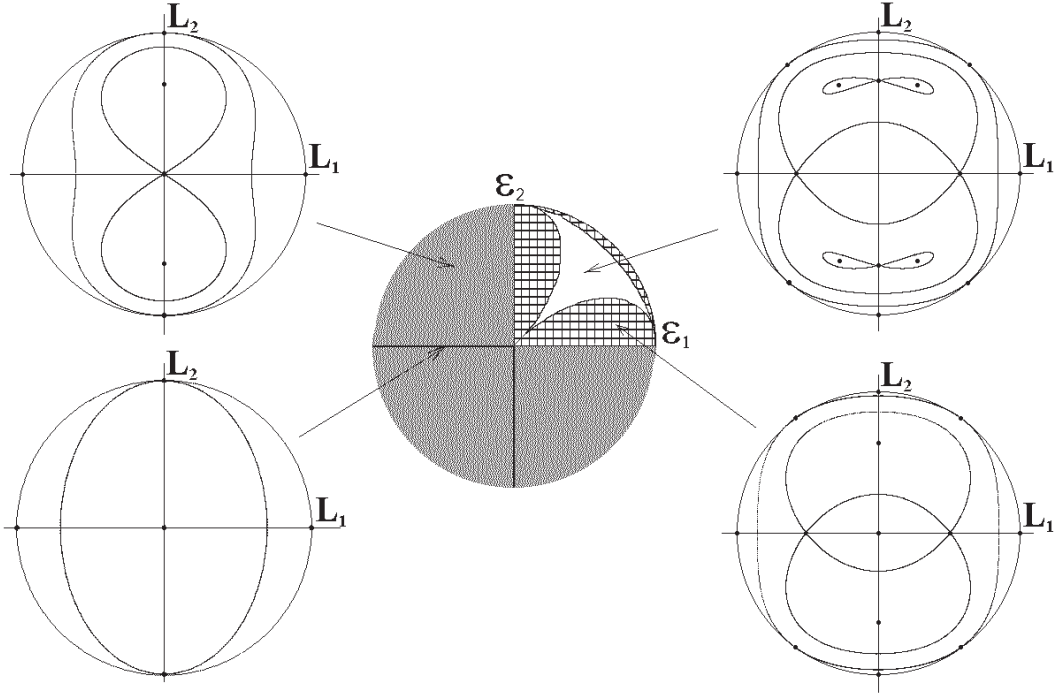}
    \caption{
         Regions of $\varepsilon_i$ corresponding to Types I - IV of
         the solutions to the auxiliary system, indicated on
         the projective plain represented by a disk with identified
         antipodal points at the boundary. The white area
         indicates {\sf Type I} solutions, the filled and the
         barred ones {\sf Type II} and {\sf Type III}.
	 The lines dividing the {\sf Type I} and {\sf Type II} regions
         are subject to Eq.(\ref{boundary})
        }
        \label{fig11}
  \end{center}
\end{figure}

\begin{figure}
  \begin{center}
    \includegraphics[width = 5in, height = 3 in] {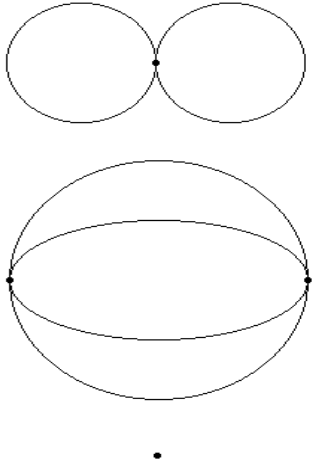}
    \caption{Topological types of the separatrixe nets;
             {\sf Type I --- IV}.
            }
            \label{fig12}
  \end{center}
\end{figure}

\begin{figure}
  \begin{center}
	\includegraphics[width = 400bp]{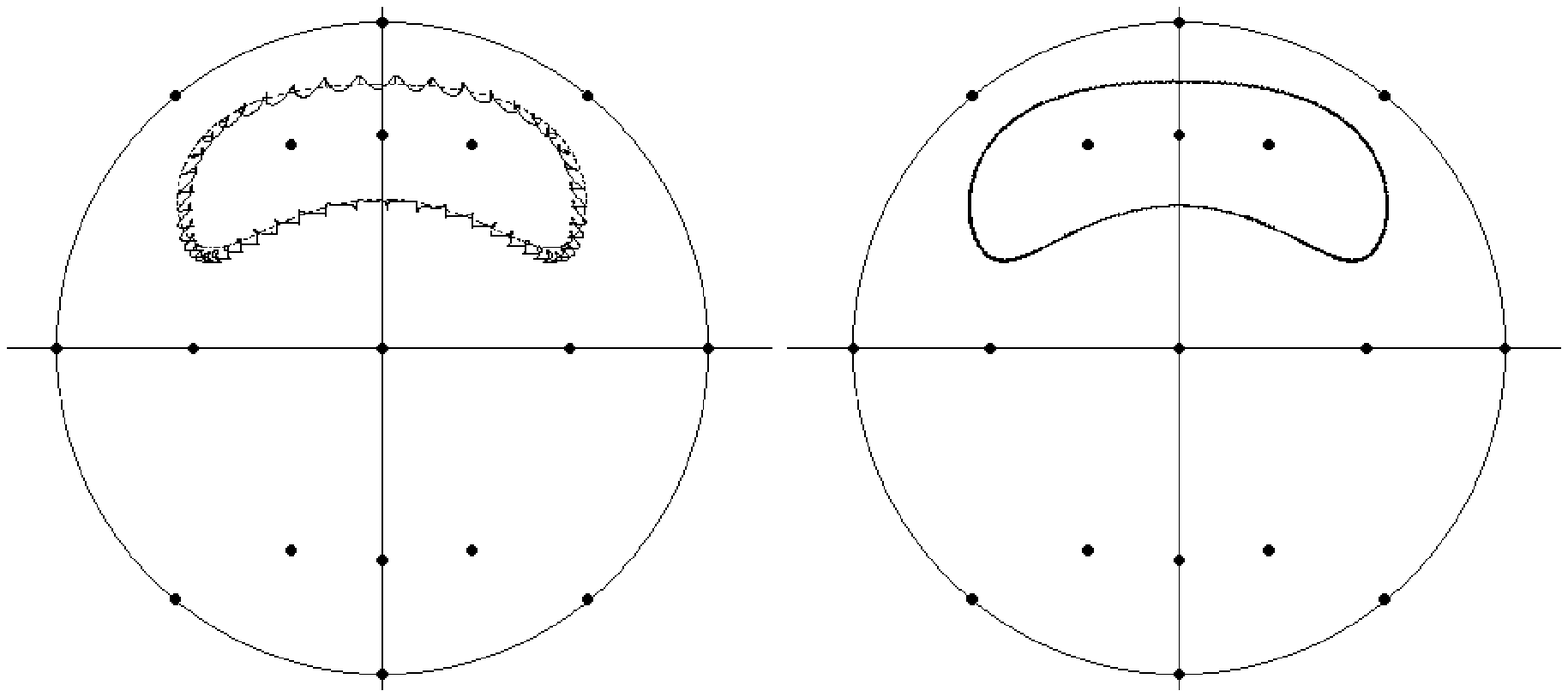}
    \caption{Comparison of the solution to: A. the initial equations
             for geodesics; B. the averaged equation given by the
             auxiliary system.
        }
        \label{fig13}
  \end{center}
\end{figure}

The topological types of the separatrix nets depend on values of
the coefficients of the deformation $\varepsilon_i$, and generate
regions I, II, III, IV in the
$\varepsilon_i$ space. Taking into account the homogeneous form of
the constraints imposed on $\varepsilon_i$, we may visualize them
on the projective plane corresponding to $\varepsilon_i$,
see FIG.\ref{fig9}.

It is important that the lines dividing the domains corresponding
to types I and II, FIG.\ref{fig11}, are given by the homogeneous equations

\begin{eqnarray}
   &1.&\quad \varepsilon_1 \varepsilon_2
  -\varepsilon_2 \varepsilon_3
  +\varepsilon_3 \varepsilon_1 = 0 \label{boundary} \\
   &2.&\quad \varepsilon_1 \varepsilon_2
  +\varepsilon_2 \varepsilon_3
  -\varepsilon_3 \varepsilon_1 = 0 \nonumber \\
   &3.&\quad -\varepsilon_1 \varepsilon_2
  +\varepsilon_2 \varepsilon_3
  +\varepsilon_3 \varepsilon_1 = 0 \nonumber
\end{eqnarray}

\begin{figure}
  \begin{center}
    \includegraphics[width = 200bp]{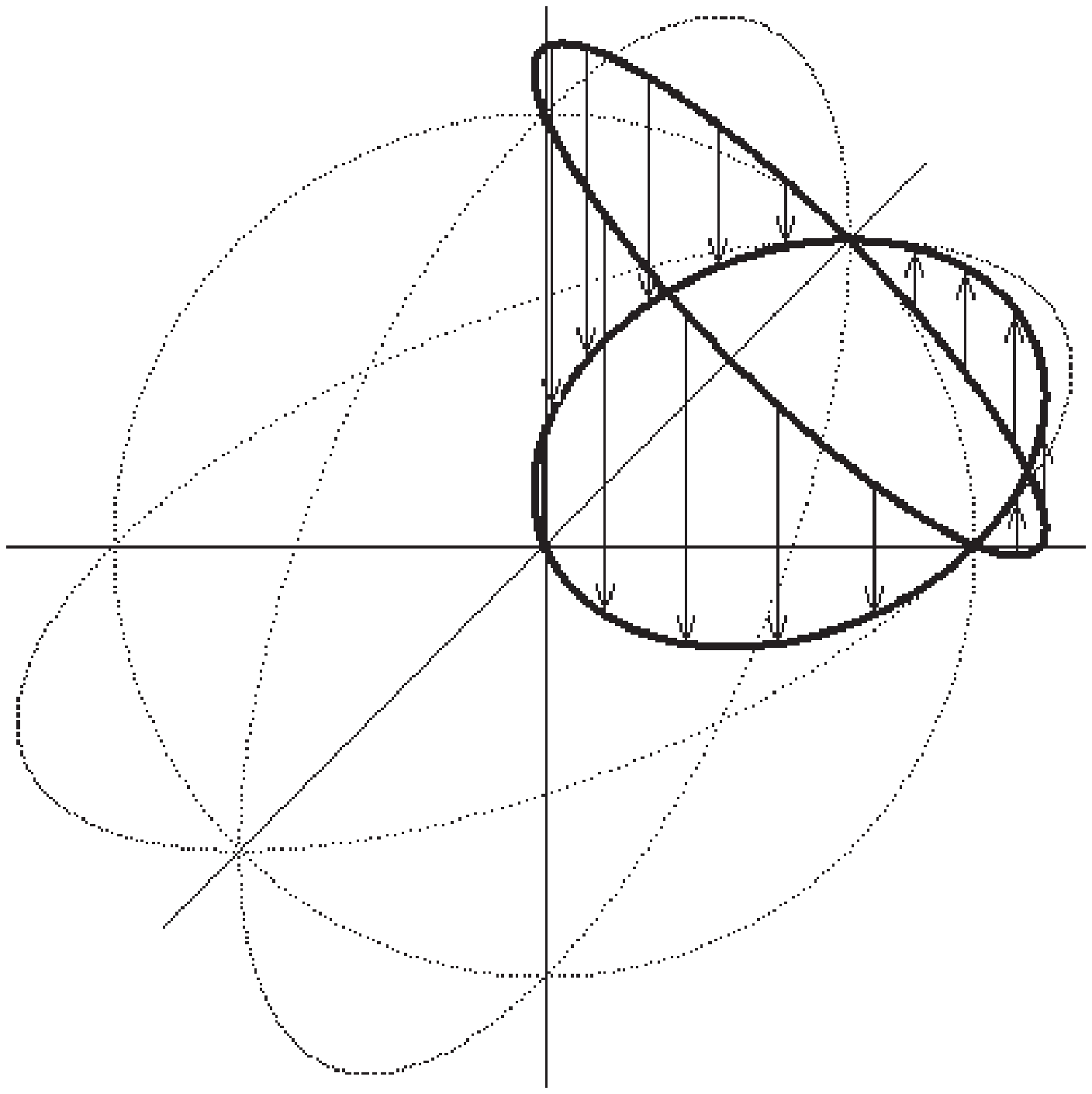}
    \caption{Intersection of the conic determining the boundary line
             of {\sf Type I} and the sphere. The projection on the disk
             gives the representation of the line on the projective plain.
            }
        \label{fig14}
  \end{center}
\end{figure}

Each equation determines a projective circle on the projective
plane of the homogeneous variables $ \varepsilon_1, \varepsilon_2,
\varepsilon_3 $, see FIGs.\ref{fig14}.

It is important that the solutions to Eqs.(\ref{avmom}) have the
specific feature that  the topological type of the graph is
completely determined by the numbers of foci and saddle points.
The dependence of the conformations of the foci and the saddles on
values of  $\varepsilon_i$ is illustrated in FIG.\ref{fig9}.

We see that the topological types of the graphs are non-trivial enough,
and the graphs realized on the sphere differ from
those on the projective plain. This circumstance is due to the fact that
the graphs on the projective plain are obtained from those on the sphere
by factoring with transformation (\ref{sym2}), see FIG.\ref{fig13}.

   \section{Conclusion}
        \label{conclusion}

The key point of the present investigation is the concept of
geodetic coil which  enables us to cast intuitive geometric ideas
in analytical form, and  relies on constructing an auxiliary
hamiltonian system, which can be considered as a reduction of the
initial problem to that of constructing a graph on projective
plane. In analytical terms, one may consider it as an asymptotic
reduction of the system of equations for geodesics on a deformed
sphere to the system similar to that of the top, with the
hamiltonian of the fourth order. The simplification we get in this
way, is substantial. Indeed, the Hamiltonian system for geodesics
could be non-integrable,whereas the auxiliary system is totally
integrable, and  its phase space can be described by a graph that
comprises vertices, which correspond to stationary solutions, or
almost closed geodesics, and edges, which can be visualized as
geodetic coils joining them. The important thing is that the
arguments, of purely analytical and topological nature, which this
analysis involves, turn out to be helpful for the numerical
simulation of the equations for geodetic lines, which admits us to
give a tangible realization of the visualization problem for
geodesics on a surface.

In fact, this paper is profoundly motivated by the technical means
provided by numerical modelling, which has allowed us to obtain
the final picture of the problem's phase space in terms of the
separatrixe graphs.  Thus, we feel that the approach used in this
paper is the symbiosis of the methods of analytical mechanics,
computer analysis, and topology. The latter is particularly
important, giving the conceptual structure for the problem of
geodetic lines on the deformed sphere. At this point it should be
noted that the detailed picture of the phase space reconstruction
should provide a description of imbedding the invariant torii,
which could  exist in some regimes, and, perhaps, the domains of
chaotic dynamics, peculiar for non-integrable problems, \cite{ab}.

The solution obtained in this paper for the geodesics on the
deformed sphere does not preclude the chaotic regimes of geodetic
lines, even though the auxiliary hamiltonian system turns out to
be completely integrable. It is worth noting that the dynamics of
geodesics gives a graphic example of the  chaos in riemanniann
geometry, \cite{ab}. This set of important problems, lying at the
border between nonlinear mechanics, and riemanniann geometry and
topology, deserves further studying.

\noindent {\bf Acknowledgment} \\
The authors are thankful to A.T.Fomenko, A.V.Bolsinov, and Yu.S.Volkov
for the useful discussions.

\noindent This work was supported by the grants NS - 1988.2003.1,
and RFFI 01-01-00583, 03-02-16173, 04-04-49645.

\end{document}